\def\versiondate{13 May 2013}
\input math.macros

\let\nobibtex = t
\let\noarrow = t
\input eplain.tex
\beginpackages
\usepackage{url}
   \usepackage{color}  %  This is for the hyperlinks.
\endpackages
\enablehyperlinks[idxexact]

\newif\ifdvipdfm   % this should be true iff getting output via dvipdfm
\input lrlEPSfig.macros
%\input Ref.macros

%\proofmodetrue
%\leftsectionheadtrue
\checkdefinedreferencetrue
%\continuousnumberingtrue
\continuousfigurenumberingtrue
\theoremcountingtrue
\sectionnumberstrue
%\figuresectionnumberstrue
\forwardreferencetrue
%\lefteqnumberstrue
%\tocgenerationtrue
\citationgenerationtrue
\nobracketcittrue
\hyperstrue
%\dvipdfmtrue
\initialeqmacro

\input\jobname.key
\bibsty{myapalike}

\def\verts{{\ss V}}
\def\vertex{{\ss V}}
\def\Td{{\Bbb T}_d}  %% regular tree of degree d
\def\Tdo{{\bf o}}  %% root of regular tree of degree d
\def\T{{\Bbb T}}  %% regular tree 
\def\F{{\scr N}}  %% networks that will be used as events
\def\Fr{{\Bbb F}}  %% free group
\def\gp{\Gamma}
\def\rtd{\mu}
\def\gh{G}
\def\edge{{\ss E}}
\def\marks{\Xi}
\def\bp{o}
\def\GG{{\scr G}}
\def\gtwo{{{\scr G}_{**}}}  % networks with 2 distinguished vertices
  % unimodular probability measures on $\GG_*$
\def\cd{\Rightarrow}  % convergence in distribution
\def\HH{{\Bbb H}}
\def\dual#1{{#1^\dagger}}
\font\frak=eufm10   %% or  eufb10
\font\scriptfrak=eufm7
\font\scriptscriptfrak=eufm5
\def\mathfrak#1#2{%       %% This cannot be done as
\def#1{{\mathchoice%
{{\hbox{\frak #2}}}%
{{\hbox{\frak #2}}}%
{{\hbox{\scriptfrak #2}}}%
{{\hbox{\scriptscriptfrak #2}}}}}}
\mathfrak{\fo}{F}  %% forest
\def\sd#1{#1^\times}     %% semi-dual object
\def\un#1#2{\underbrace{#1}_{\rm (#2)}}

\ifproofmode \relax \else\head{To appear in {\it Ergodic Theory Dynamical
Systems}}
{Version of \versiondate}\fi 
\vglue20pt

\title{Unimodular Random Trees}

\author{Itai Benjamini, Russell Lyons, and Oded
Schramm{\eightpoint\parindent=10pt\footnote{${}^\dagger$}{Deceased.}}}

\abstract{We consider unimodular random
rooted trees (URTs) and invariant forests in Cayley graphs.
We show that URTs
of bounded degree are the same as the law of the component of the root in
an invariant percolation on a regular tree. We use this to give a new proof
that URTs are sofic, a result of Elek. We show that ends of
invariant forests in the hyperbolic plane converge to ideal boundary points.
We also note that uniform integrability of the degree distribution of a
family of finite graphs implies tightness of that family for local
convergence, also known as random weak convergence.
}

\bottomIII{Primary 
 05C05, %   Trees
60C05. % Combinatorial probability
Secondary
 05C80, %   Random graphs
82B43, % Percolation
60G10. %Stationary processes
}
{Graphs, invariant percolation, sofic, treeable, groups, ends, forests,
hyperbolic space, random weak convergence, tight, degree distribution.}
{Research partially supported by Microsoft Research and
NSF grant DMS-1007244.}

\vskip-10pt

\bsection{Introduction}{s.intro}

The theory of unimodular random rooted networks (URNs) is an outgrowth
mainly of two lines of investigations: one is concerned with
asymptotics of finite networks, while the other involves group-invariant
stochastic processes on infinite Cayley
graphs, especially percolation.
An important motivation also arises from the class of
sofic groups. Parallels with the theory of limits of dense graphs now spur
further investigations (see \ref b.Lovasz:limits/ for this).

We give full definitions in \ref s.def/, but here we recount
intuitively some of the above motivations. 
One way to look at a large finite network (which is a labeled graph) is to
look at a large neighborhood around a random uniformly chosen vertex.
Often such neighborhood statistics capture quantities of interest and their
asymptotics.
Thus, one is led to take limits of such statistics and thereby define a
probability measure on infinite rooted graphs, where the neighborhood of
the root has the statistics that arise as the limit statistics of the
finite networks.
Such a limit of a sequence of finite networks is called the random weak
limit, the local (weak) limit, the distributional limit, or the
Benjamini-Schramm limit of the sequence.
All such limit measures have a property known as unimodularity; it is not
known whether all unimodular measures are limits of finite networks.
Those that are such limits are called sofic.

Intuitively, a probability measure on rooted networks is unimodular iff its
root is chosen ``uniformly" from among all its vertices.
This, of course, only makes sense for finite graphs. It is formalized for
networks on infinite graphs by requiring a sort of conservation property
known as the Mass-Transport Principle.

Unimodularity is an extremely powerful property,
especially for studying percolation on infinite graphs. In the present
context,
for example, the component of the identity in a group-invariant
percolation on a Cayley graph has a unimodular law as a random
graph rooted at the identity.

Consider the following example of a random weak limit of finite graphs: Let
$\T_3$ be a 3-regular tree. Let $G_n$ be the ball of radius $n$ in $\T_3$
about any point. Since most points in $G_n$ are near the leaves of $G_n$,
the random weak limit of $\Seq{G_n \st n \ge 1}$ is not $\T_3$ but the
following probability measure, $\mu$. 
Let $T_n$ be disjoint binary trees of depth $n$ for $n \ge 0$.
Modify $T_n$ by adding a new vertex $x_n$ adjacent
to the root of $T_n$. Also consider
an isolated vertex $x_{-1}$. Now add an edge between
$x_n$ and $x_{n+1}$ for each $n \ge -1$. The resulting tree,
$T$, was called the \dfn{canopy tree} by \ref b.AziWar/. 
The graph $T$ rooted at $x_n$ is denoted $(T, x_n)$.
We now define $\mu$ by letting $\mu(T, x_{n}) := 2^{-n-2}$ ($n \ge -1$).
Thus, $\mu$ is supported on a single tree, which is a proper subtree of
$\T_3$.
In fact, as we show, $\mu$ can be obtained as the component
of the root in an automorphism-invariant percolation on $\T_3$.
Indeed, one of our main theorems is that every URN that is supported by
trees of bounded degree can be obtained as the component of the root in an
invariant percolation on a regular tree.

An interesting contrast is provided by other URNs.
For example, consider the infinite discrete Sierpi\'nski gaskets
characterized by \ref b.Tepl/ (see Lemma 2.3 there).
These are obtained as the random weak limit of the graphs that are the
natural boundaries of the $n$th-stage construction
of the usual Sierpi\'nski gasket as $n \to\infty$ (see \ref
f.gasket-graphs/).
In this case, the limit measure $\mu$ has an uncountable support, although
still all the graphs in the support are subgraphs of the triangular lattice in
the plane.
Yet in this case,
there is no invariant percolation on the triangular lattice such that
the component of the root has law $\mu$ since the subgraphs in the support
of $\mu$ have density 0 and by their topology, only one component in any
percolation on the triangular lattice can be a Sierpi\'nski gasket (except
for degenerate ones that altogether have $\mu$-measure 0).
%Although the Menger sponge limits don't have this same property relative
%to $\Z^3$, they too cannot be obtained by invariant percolation, as can be
%seen by considering 2D slices.

\efiginslabel gasket-graphs {The first three graphs whose random weak limit is the
infinite Sierpi\'nski gasket.} y1.5

We call a URN that is supported by trees a URT.
We shall use our result that every URT of bounded degree can be obtained as
the law of the component of the root in an invariant labeled percolation on
a regular tree to give a new proof that URTs are sofic.
This was first shown (in a special case) by
\ref b.Elek:sofictrees/, which solved
Question 3.3 of \ref b.BolRio/. It was then extended from graphs to networks
by \ref b.ElekLipp:SER/.
Although not needed for any of these results, we give a sufficient
condition for a collection of finite graphs to have a convergent
subsequence, namely, that their degree distributions be uniformly
integrable.

We remark that our theorem showing that every URT of bounded degree can be
obtained via invariant percolation on a regular tree has a counterpart in
the other direction: That is, rather than put a URT into a regular tree,
one can put a regular tree (or forest) on a URT. More precisely, \ref
b.Hjorth/ proved that every treeable (probability-measure-preserving)
equivalence relation of cost at least 2 can be generated by a free action
of a free group $\Fr_2$ on 2 generators. If, instead, the cost is assumed
only to be larger than 1, then there is a subrelation that is generated by a
free action of $\Fr_2$: see Proposition 14 of \ref b.GabLy/. In the remaining
case
where the treeable
equivalence relation has cost 1, the equivalence relation is amenable,
whence a theorem of \ref b.CFW/ shows that it is generated by a
free action of $\Z$. A URT is essentially the same as a
treeable equivalence relation.  We do not use any of these notions here,
however, so we leave these terms undefined.

Finally, we turn from results about general URTs to somewhat specific URTs.
Consider a discrete forest in the hyperbolic plane.
A simple infinite path in the forest is called a ray. Clearly a ray can have
fairly arbitrary limiting behavior; in particular, though it must tend to
the ideal boundary because the forest is discrete, the ray need not converge
to any ideal boundary point.
However, we show that with the condition solely that the forest is random
with a law that is invariant under hyperbolic isometries, a.s.\ all its rays
converge to ideal boundary points.
We do not know whether this holds in higher dimensions or, more generally,
in word-hyperbolic groups.
We do know that rays do not necessarily converge with positive speed.
Note that if we fixed a point in the hyperbolic plane and took the nearest
vertex of the forest to that point as the root of its component, then we
would obtain a URT.

\bsection{Definitions}{s.def}

We review a few definitions from the theory of unimodular random rooted
networks; for more details, see \ref b.AL:urn/.
A \dfn{network} is a (multi-)graph $\gh = (\vertex, \edge)$ together with a
complete separable metric space $\marks$ called the \dfn{mark space} and
maps from $\vertex$ and $\edge$ to $\marks$. Images in $\marks$ are called
\dfn{marks}.
Each edge is given two marks, one associated to (``at") each of its endpoints.
The only assumption on degrees is that they are finite.
We omit the mark maps from our notation for
networks.

A \dfn{rooted network} $(\gh, \bp)$ is a network $\gh$ with a distinguished
vertex $\bp$ of $\gh$, called the \dfn{root}.
A \dfn{rooted isomorphism} of rooted networks is an isomorphism of the
underlying networks that takes the root of one to the root of the other.
We do not distinguish between a rooted network and its
isomorphism class.
Let $\GG_*$ denote the set of rooted isomorphism classes of rooted
{\it connected\/} locally finite networks.
Define a separable complete metric on $\GG_*$ by letting the distance
between $(G_1, o_1)$ and $(G_2, o_2)$ be $1/(1+\alpha)$, where $\alpha$ is
the supremum of those $r > 0$ such that there is some rooted isomorphism of
the balls of (graph-distance) radius $\flr{r}$ around the roots of $G_i$
such that each pair of corresponding marks has distance less than $1/r$.
For probability measures $\rtd$, $\rtd_n$ on $\GG_*$, we write $\rtd_n \cd
\rtd$ when $\rtd_n$ converges weakly with respect to this metric.

For a (possibly disconnected)
network $\gh$ and a vertex $x \in \verts(\gh)$, write $\gh_x$ for the
connected component of $x$ in $\gh$.
If $\gh$ is finite, then write $U_\gh$ for a uniform random vertex of $\gh$
and $U(\gh)$ for the
corresponding distribution of $\big(\gh_{U_\gh}, U_\gh\big)$ on $\GG_*$.
Suppose that $\gh_n$ are finite networks and that $\rtd$ is a
probability measure on $\GG_*$.
We say that the \dfn{random weak limit} of $\gh_n$ is $\rtd$ if $U(\gh_n)
\cd \rtd$.

A probability measure
that is a random weak limit of finite networks is called \dfn{sofic}.
In particular, a finitely generated
group is called sofic when its Cayley diagram is sofic.
It is easy to check that this property does not depend on the generating
set chosen.
All sofic measures are unimodular, which we now define.
Similarly to the space $\GG_*$, we define the space $\gtwo$ of isomorphism
classes of locally
finite connected networks with an ordered pair of distinguished vertices
and the natural topology thereon:
the distance
between $(G_1, o_1, o'_1)$ and $(G_2, o_2, o'_2)$ is
$1/(1+\alpha)$, where $\alpha$ is
the supremum of those $r > 0$ such that there is some isomorphism of
the balls of radius $\flr{r}$ around $o_i$ that takes
$o_1$ to $o_2$ and
$o'_1$ to $o'_2$
such that each pair of corresponding marks has distance less than $1/r$.
We shall write a function $f$ on $\gtwo$ as $f(\gh, x, y)$.
We refer to $f(\gh, x, y)$ as the \dfn{mass} sent from $x$ to $y$ in $\gh$.

\procl d.unimodular 
Let $\rtd$ be a probability measure on $\GG_*$.
We call $\rtd$ \dfn{unimodular} if it obeys the \dfn{Mass-Transport
Principle}:
For all Borel
$f : \gtwo \to [0, \infty]$,
we have
$$
\int \sum_{x \in \vertex(\gh)} f(\gh, \bp, x) \,d\rtd(\gh, \bp)
=
\int \sum_{x \in \vertex(\gh)} f(\gh, x, \bp) \,d\rtd(\gh, \bp)
\,.
\label e.mtpgen
$$
\endprocl

It is easy to see that every sofic measure 
is unimodular, as observed by \ref b.BS:rdl/, who introduced this
general form of the Mass-Transport Principle under the name 
``intrinsic Mass-Transport Principle".
The converse is open and was posed as a question by \ref b.AL:urn/.

A special form of the Mass-Transport Principle was considered, in different
language, by \ref b.AS:obj/.
Namely, they defined $\rtd$ to be \dfn{involution invariant} if \ref
e.mtpgen/ holds for those $f$ supported on $(\gh, x, y)$ with $x \sim y$.
In fact, the Mass-Transport Principle holds for general $f$ if it holds for
these special $f$, as shown by \ref b.AL:urn/:

\procl p.invinv
A measure is involution invariant iff it is unimodular.
\endprocl

If $\gh$ is $\rtd$-a.s.\ regular, then involution invariance of $\rtd$ is
equivalent to the following: If $o'$ is a uniform random neighbor of the
root, then the law of $(\gh, o, o')$ is the same as the law of $(\gh, o',
o)$ when $(\gh, \bp)$ has the law $\rtd$.

See also \ref b.BC:stat/ for a discussion of unimodularity.

We call a measure a \dfn{URT} if it is a unimodular probability measure on
rooted networks whose underlying graphs are trees.
We call a probability measure a \dfn{labeled percolation} on a graph $G$ if
it is carried by the set of networks on $G$ whose marks are pairs, with the
second coordinate, called \dfn{color},
of a mark being 0 or 1.  Edges colored 0 or 1 are
called \dfn{closed} and \dfn{open}, respectively.

\bsection{Tightness and Degree}{s.tight}

One of our theorems is that URTs are sofic. Although for this purpose, we
shall not need results that imply random weak convergence of a subsequence
of finite graphs, such results have not been stated in the literature before
except in the easy case of bounded degree and the harder case of
exponential tails of the degree distribution (\ref B.AS:uipt/).
On the other hand,
it does not suffice that the mean degrees be bounded. For example, consider
the complete bipartite graphs $K_{1, n}$ (stars):
No subsequence converges.
Yet also, it is not necessary that the mean degrees be bounded.
In fact, the mean degree can be infinite for extremal sofic unimodular
random rooted graphs, even trees.
Here, \dfn{extremal} means that the probability measure is not a convex
combination of other unimodular probability measures on rooted graphs.
We impose that condition since it is trivial to take a mixture
of finite-mean-degree URTs to get a URT of infinite mean degree.

For examples, let $\Seq{p_n \st n \ge 1}$ be a probability distribution on
$\Z^+$ with infinite mean.
For each integer $k$, join $k$ to $k + 1$ by $n$ parallel edges with
probability $p_n$, independently for different $k$.
This is easily seen to be an extremal sofic probability measure.
To get an extremal URT with infinite mean degree, take
the universal cover rooted at 0 of the resulting multigraph;
see Example 9.3 of \ref b.AL:urn/.

% One could also put a K_n at k with probability p_n, rooted at a uniform
% vertex of the clique at 0.
Thus, it may be useful to present the following result on tightness.
For simplicity and with no essential loss of generality,
we shall assume that all our graphs have no isolated vertices.
The proof we present was suggested by Omer Angel and simplifies our
original proof. 
The proof in \ref b.AS:uipt/ would also work.

\procl t.tight
If $A$ is a family of finite graphs such that the random variables
$\big\{\!\deg_G U_G \st G \in A\big\}$ are uniformly integrable, then
$\big\{U(G) \st G \in A\big\}$ is tight.
\endprocl

\proof
Let $f(d) := \sup_{G \in A} \E\big[\!\deg_G U_G \st \deg_G U_G > d\big]$.
By assumption, $\lim_{d \to\infty} f(d) = 0$.
Write $m(G) := \E\big[\!\deg_G U_G\big]$.
Thus, $1 \le m(G) \le f(0) < \infty$.
Write $D(G)$ for the degree-biased probability measure on $\big\{(G,
x) \st x \in V(G) \big\}$, that is, 
$$
D(G)\big[(G, x)\big] = {\deg_G x \over m(G)}
\cdot U(G)\big[(G, x)\big]
\,,
$$
and $D_G$ for the corresponding root.
Since $U(G) \le m(G) D(G) \le f(0) D(G)$, it suffices to show that $\big\{D(G)
\st G \in A\big\}$ is tight.
Note that $\big\{\!\deg_G D_G \st G \in A\big\}$ is tight.

For $r \in \N$, let $F_r^M(x)$ be the event that there is some vertex at
distance at most $r$ from $x$ whose degree is larger than $M$.
Let $X$ be a uniform random neighbor of $D_G$.
Because $D(G)$ is a stationary measure for simple random walk, 
$F_r^M(D_G)$ and $F_r^M(X)$ have the same probability.
Also, $\P[F_{r+1}^M(D_G) \mid \deg_G D_G] \le (\deg_G D_G) \P\big[F_r^M(X) \mid
\deg_G D_G\big]$.
We claim that for all $r \in \N$ and $\epsilon > 0$, there exists $M <
\infty$ such that $\P\big[F_r^M(D_G) < \epsilon\big]$ for all $G \in A$;
this clearly implies that $\big\{D(G) \st G \in A\big\}$ is tight.
We prove the claim by induction on $r$.

The statement for $r = 0$ is trivial.
Given that the property holds for $r$, let us now show it for $r+1$.
Given $\epsilon > 0$, choose $d$ so large that $\P[\deg_G D_G > d] <
\epsilon/2$ for all $G \in A$.
Also, choose $M$ so large that $\P\big[F_r^M(D_G)\big] < \epsilon/(2d)$
for all $G \in A$.
Write $F$ for the event that $\deg_G D_G > d$.
Then by conditioning on $\deg_G D_G$, we see that
$$\eqaln{
\P[F_{r+1}^M(D_G)]
&\le
\P[F] + \E\big[\I{F^{\rm c}} \P[F_{r+1}^M(D_G) \mid \deg_G
D_G]\big]
\cr&\le
\epsilon/2 + \E\big[\I{F^{\rm c}} (\deg_G D_G) \P[F_r^M(X) \mid \deg_G
D_G]\big]
\cr&\le
\epsilon/2 + \E\big[\I{F^{\rm c}} d \P[F_r^M(X) \mid \deg_G D_G]\big]
\cr&\le
\epsilon/2 + \E\big[d \P[F_r^M(X) \mid \deg_G D_G]\big]
\cr&=
\epsilon/2 + d \P\big[F_r^M(X)\big]
\cr&=
\epsilon/2 + d \P\big[F_r^M(D_G)\big]
\cr&<
\epsilon/2 + d\epsilon/(2d)
= \epsilon
}$$
for all $G \in A$, which proves the claim.
\Qed

In this proof, the only way that we used finiteness of the graphs was that
the degree-biased uniform distribution on vertices gave a stationary
measure for simple random walk. Thus, the result applies also to any
collection of
probability measures bounded by a fixed multiple of
stationary probability measures on rooted graphs, such as unimodular
probability measures on graphs.

\bsection{Invariant Percolation}{s.invar}

We now prove that every URT of bounded degree arises as the open component
of the root in an invariant percolation on a regular tree.

We use the following lemma that
is straightforward to check from the definitions.
The technical definition of adding IID marks is explained in Section 6 of
\ref b.AL:urn/.

\procl l.unimodular-factor
Suppose that $\mu$ is a unimodular probability measure on rooted networks.
Let $\phi$ be a measurable map on rooted networks that takes each network to
an element of the mark space. Define $\Phi$ to be the map on rooted
networks that takes a network $(G, o)$ to another network
on the same underlying graph, but replaces the mark at each vertex $x \in
G$ by $\phi(G, x)$.
Then the pushforward measure $\Phi_*\mu$ is also unimodular.
If instead we add a second coordinate to each vertex mark by an IID mark
according to some probability measure on the mark space,
then the resulting measure is again unimodular.
\endprocl

\procl t.urt-invar
Let $\mu$ be a probability measure on rooted networks whose underlying
graphs are trees of degree at most $d$. Then $\mu$ is unimodular iff $\mu$
is the law of the open component of the root in a labeled percolation on a
$d$-regular tree whose law is invariant under all automorphisms of the tree.
\endprocl

\proof
The ``if" part of the assertion is well known and not dependent on the fact
that the underlying graph is a tree.
See, e.g., \ref b.BLPS:gip/ or Theorem 3.2 of \ref b.AL:urn/.

The idea for proving the converse is as follows.
First sample $(T, o) \sim \mu$. Of all the possible ways to embed it in the
$d$-regular tree $\Td$ such that $o$ maps to the root $\Tdo$ of $\Td$, choose
one uniformly (i.e., choose one arbitrarily and then apply a uniform
automorphism of $\Td$ preserving $\Tdo$). The embedded
image of $T$ is marked open.
Now for every
edge $e$ in $\Td$ that is not in the image of
$T$ but has one endpoint in $T$, mark $e$ closed and
sample an independent copy $(T', o') \sim \mu$ with 
$o'$ embedding as the endpoint of $e$ that is not in $T$. However,
this choice of $T'$ has to be biased so that the degree of $o'$
is not $d$. In fact, we sample instead $(T', o') \sim \mu'$, where $\mu'$ is
absolutely continuous with respect to $\mu$ with Radon-Nikod\'ym derivative
at $(T, o)$ equal to $(d - \deg_T o)/\alpha$, where $\alpha$ is a normalizing
constant.
Continue in this way to cover all of
the vertices of $\Td$ by weighted independent
copies of $(T', o') \sim \mu'$. 
Of course, all edges in embedded copies of $T$ or $T'$ are marked open,
while the rest are marked closed.
To prove that this is invariant, 
we first show involution invariance of the constructed marked tree and then
appeal to Theorem 3.2 of \ref b.AL:urn/
to get that it is an invariant percolation on $\Td$.
Proving involution invariance involves two cases: one case involves
crossing a closed edge; that's where the biased measure $\mu'$ comes in.
The other case
involves crossing an open edge; that's where the unimodularity of the
original measure $\mu$ comes in.

Here are the details.
Let $p_k$ be the $\mu$-probability that the root has $k$ children. 
Let $\alpha := \sum_k p_k (d-k)$ be a normalizing constant.
Let $\mu'$ be
absolutely continuous with respect to $\mu$ with Radon-Nikod\'ym derivative
at $(T, o)$ equal to $(d - \deg_T o)/\alpha$.
Given two probability measures $\nu_1$ and $\nu_2$ supported by networks on
rooted trees where the root has degree at most $d-1$ and all other vertices
have degree at most $d$, write $Q(\nu_1, \nu_2)$ for the probability
measure supported by networks on the rooted $d$-ary tree
constructed as follows, similar to a Galton-Watson branching process:
Choose $(T', o) \sim \nu_1$, whose edges are colored open.
To each vertex $x \ne o$ of $T'$, adjoin $d - \deg_{T'} x$ edges colored
closed
whose other endpoint is the root of an independent sample from $\nu_2$, while
to the root $o$ of $T'$, adjoin $d - 1 - \deg_{T'} o$ edges colored
closed 
whose other endpoint is the root of an independent sample from $\nu_2$.
Call the result the network $(T, o)$. Then $Q(\nu_1, \nu_2)$ is the law of
$(T, o)$.
Write $\nu$ for the measure on rooted networks defined by the equation
$\nu = Q(\mu', \nu)$.

Let
$\rho$ be the measure constructed as follows:
Choose $(T', o) \sim \mu$, whose edges are colored open.
To each vertex $x$ of $T'$, adjoin $d - \deg_{T'} x$ edges colored
closed
whose other endpoint is the root of an independent sample from $\nu$.
The result is a network whose underlying graph is $\Td$.
We claim that this measure $\rho$ is unimodular, which we show by
proving that $\rho$ is involution invariant.
This suffices to prove the theorem by appeal to Theorem 3.2 of \ref
b.AL:urn/ (in which the averaging over automorphisms is taken).

To prove this claim, 
it will be convenient to use the following technical modification of $\rho$
to deal with counting issues:
% Note that this concerns not only permutations of the $B_m$ below, but
% also permutations that involve $(o, o')$, which is quite messy.
Given $(T, o) \sim \rho$, assign independently and uniformly marks to the
closed edges in each direction so that each vertex is surrounded by
outgoing closed edges marked $1, \ldots, k$ when it is incident to $k$
closed edges.
Call $\rho'$ the resulting measure on networks.
It clearly suffices to prove that $\rho'$ is involution invariant.

For $k$ such that $p_k > 0$, let $\mu_k$ be the measure constructed as follows:
Choose $(T', o) \sim \mu$, whose edges are colored open, conditioned on
$\deg_{T'} o = k$.
To each vertex $x \ne o$ of $T'$, adjoin $d - \deg_{T'} x$ edges colored
closed
whose other endpoint is the root of an independent sample from $\nu$.
Let $\F_i$ denote the class of
networks supported on a rooted tree with all vertices
except the root having degree $d$ and the root having degree $i$.
Consider $i, i' \in [1, d-1]$ and Borel sets $A \subseteq \F_i$, $A' \subseteq
\F_{i'}$, and $B_1, \ldots, B_{d-i-1}, B'_1, \ldots, B'_{d-i'-1} \subseteq
\F_{d-1}$.
Now let $(T, o) \sim \rho'$ and let $o'$ be a uniform neighbor of $o$.
Then the chance that we see (a) $i$ open edges at $o$, (b)
the edge $(o, o')$ is
closed with (c) mark $j$ in the direction $(o, o')$ and (d) 
mark $j'$ in direction
$(o', o)$, see (e) $i'$ open edges at $o'$, and see 
the event where (f) the open
edges at $o$ are part of a network in $A$, (g) the open edges at $o'$ are part
of a network in $A'$, while (h) the other endpoints of the closed edges at $o$
belong to networks in $(B_1, \ldots, B_{d-i-1})$ in increasing order of
their marks from $o$ and (i) similarly the other endpoints of the closed edges
at $o'$ belong to networks in $(B'_1, \ldots, B'_{d-i'-1})$ in increasing
order of their marks from $o'$
(see \ref f.complexevent/)
equals
$$
\un{p_i}{a} \cdot \un{\mu_i(A)}{f} \cdot \un{\prod_{m=1}^{d-i-1} \nu(B_m)}{h}
\cdot \un{d-i \over d}{b} \cdot
\un{1 \over (d - i) }{c} \cdot
\un{1 \over  (d - i')}{d}
\cdot \un{\prod_{r=1}^{d-i'-1}
\nu(B'_r)}{i} \cdot \un{p_{i'} {d-i' \over \alpha}}{e} \cdot
\un{\mu_{i'}(A')}{g}\,.
$$
This is invariant under the involution exchanging $o$ and $o'$.

\lrlfiginslabellong complexevent {All edges drawn with solid lines
are closed. The root $o$ is
incident to 3 closed edges, while $o'$ is incident to 4 closed edges.} 4.5 x4
{
%\ShowGrid
(.585*.38) $o$\\
(.42*.38) $o'$\\
\E(.85*.25) $A$\\
\E(.15*.35) $A'$\\
\E(.9*.6) $B_1$\\
\L(.9*.8) $B_2$\\
\E(.1*.55) $B'_1$\\
\E(.1*.72) $B'_2$\\
\E(.1*.9) $B'_3$\\
}

The other case to prove is when the edge $(o, o')$ is open.
Consider the following measure, $\sigma$. 
Begin with a sample $(T, o) \sim \mu$.
Assign a second coordinate $\big(\tau_1(x), \ldots, \tau_{d-1}(x)\big)$
to the vertex mark at each vertex $x$ given by IID samples $\tau_i(x) \sim
\nu$.
This new network is unimodular by \ref l.unimodular-factor/.
Now replace the second coordinate of the vertex mark at each vertex $x$ by
$\big(\tau_1(x), \ldots, \tau_{d -\deg_T(x)}(x)\big)$.
This new network is again unimodular by \ref l.unimodular-factor/.
We denote by $\sigma$ its law.
Note that we can obtain $\rho$ from $\sigma$ by replacing the second
coordinate of the vertex mark at each $x$ by a tree network rooted at $x$,
where we adjoin $d - \deg_T(x)$ closed edges to $x$, at the other end of
which we adjoin the trees $\tau_i(x)$.

What remains to prove is that involution invariance holds for $\rho'$ across
open edges. It suffices to do the same for $\rho$.
But this is clearly the same as unimodularity of $\sigma$. 
\Qed

\bsection{Soficity}{s.sofic}

We now use the preceding theorem to prove that URTs are sofic.
This result is the same as Theorem 4 of 
\ref b.ElekLipp:SER/, but in different language.
See Example 9.9 of \ref b.AL:urn/ for a comparison of the different
languages.\ftnote{${}^*$}
{Actually, the notion of ``sofic" in \ref b.ElekLipp:SER/ applies only to
certain URNs, namely, those where the rooted network is a measurable
function of the label of the root. Thus, the result we prove is
superficially more general.}

\procl t.URT-sofic
Every URT is sofic.
\endprocl

\proof
It follows from Theorem 3.4 of \ref b.Bowen:per/ that every invariant
network on $\Td$ is sofic; the result is stated there for even $d$ only,
but the proof works for all $d$.
Thus, given a URT $\mu$, if the degrees are bounded by $d$, let $\rho$ be
an invariant labeled percolation on $\Td$ such that the open component of
the root has law $\mu$.
Let $\Seq{G_n \st n \ge 0}$ be finite networks whose random weak limit is
$\rho$.
Here, we may assume that the edges of $G_n$ are each colored closed or
open.
Let $G'_n$ be the result of deleting every closed edge from $G_n$. Then
clearly $\Seq{G'_n}$ has random weak limit $\mu$.
Finally, if the degrees are not bounded $\mu$-a.s., then for each $d$, let
$\mu_d$ be the law of the component of $o$ when
we delete every edge of $T$ incident to some vertex of degree larger than
$d$, where $(T, o) \sim \mu$.
Then $\mu_d$ is unimodular and, by what we just proved, sofic.
Since the sofic measures form a weakly closed set and $\mu_d$ tend weakly
to $\mu$, we deduce that $\mu$ is sofic as well. 
\Qed

As noted by \ref b.ElekLipp:SER/, this implies that every treeable group is
sofic.
Here, a group $\gp$ is \dfn{treeable} if there is a probability measure on
trees with vertex set $\gp$ that is invariant under the natural action of
$\gp$; such a probability measure is called a \dfn{treeing} of $\gp$.
Briefly, the idea is to use a treeing $\mu$ of $\gp$, a generating
set $S$ for $\gp$, and a sofic approximation $\Seq{G_n}$ of $\mu$ to
construct a sofic approximation of the Cayley diagram of $\gp$ with respect
to $S$ by putting
edges labeled $s \in S$ between points $x, y$ of $G_n$ such that a path
from $x$ to $y$ has length at most $R_n$ and has labels that multiply to
$s$, where $R_n \to \infty$ at an appropriately slow rate.

\bsection{Rays}{s.ends}

Random weak limits of finite trees have mean degree at most
2, are supported by trees with at most 2 ends, and hence are recurrent for
simple random walk; see Proposition 6.3 of \ref b.AL:urn/.
In the case of URTs with 
finite mean degree larger than 2, the speed of simple random
walk is positive: see Theorem 4.9 of \ref b.AL:urn/.
The case of infinite mean degree is open.
However, it is interesting in all cases to see whether the rays themselves,
rather than simple random walk,
have positive speed when embedded in a larger graph.
What we mean by this is the following.

We say that a sequence $\Seq{x_n \st n \ge 0}$ in a metric space has
\dfn{positive (liminf) speed} if there is some constant $c > 0$ such that
the distance between $x_n$ and $x_0$ is at least $c n$ for all $n \ge 1$.
A simple infinite path in a tree is called a \dfn{ray}.
An \dfn{end} of a tree is an equivalence class of rays, where two rays are
equivalent when they have finite symmetric difference.
Of course, any statement about limits of rays applies equally to all rays
belonging to the same end and is therefore a statement about limits of
ends.

We are interested in the rays in forests that arise either as invariant
percolation on a Cayley graph or as random graphs discretely
embedded in hyperbolic
space $\HH^d$ with an isometry-invariant law.
When do all the rays have positive speed in the metric of the Cayley graph
or in the hyperbolic metric? It does not
suffice that the Cayley graph be non-amenable: For example,
consider the usual Cayley graph of the group
$\Z * \Z^2$. Use the random forest that arises from
an independent copy of the uniform spanning tree (\ref B.Pemantle/) in
every copy of $\Z^2$. Then almost surely, each such tree contains only one
end and no ray has positive speed. (In fact, the $n$th vertex in a ray is
roughly at distance $n^{4/5}$ from its starting point; see \ref
b.BarMass:LERW2D/.)
%Thus, we must strengthen the assumption on the Cayley graph if we want the
%conclusion to hold.
What if we restrict ourselves to word-hyperbolic groups?
As we shall see, the answer is still no.
Thus, we focus on the following weaker property for hyperbolic groups:

\procl q.rays-hyp-space
Does every ray in an invariant forest in a word-hyperbolic group
converge a.s.\ to an ideal boundary point?
\endprocl

See the survey \ref b.KapBen/ for information on the boundary of a
word-hyperbolic group.
We know the answer only in $\HH^d$ for $d = 2$:

\procl t.rays-hyp-space
Let $G$ be a one-ended graph embedded in $\HH^2$ such that a
group of isometries of $\HH^2$ acts quasi-transitively on $G$.
Given an automorphism-invariant forest in $G$, a.s.\
every ray in the forest converges to an ideal boundary point.
Furthermore, the set of limits of the rays
is a.s.\ the entire ideal boundary.
\endprocl

Here, to say that $G$ is one-ended means that the complement of each finite
set in $G$ has only one infinite component in $G$.
Note that the ideal boundary points of $G$ are the same as those of
$\HH^2$.
We call the set of limit points of the convergent
rays in a tree or forest the \dfn{limit set} of that tree or forest.

One can prove a similar statement for forests in $\HH^2$ whose law is
invariant under isometries of $\HH^2$, without assuming an underlying
graph, $G$.
On the other hand,
one could also let the underlying graph
$G$ be random with isometry-invariant law; no
quasi-transitivity of $G$ is then needed, nor need $G$ have only one end.
The lengths of edges should be bounded and the vertices should be separated
by a minimum distance in $\HH^2$.

We are grateful to Omer Angel for some simplifications to our proof.

\proof
Without loss of generality, we may assume that the forest $\fo$ is spanning
and contains only infinite trees. 
Indeed, we may first delete all finite trees and then independently
add an edge at random from each vertex not in the forest to a vertex that
is closer to the forest.

Let $\dual G$ be the planar dual graph of $G$; its edges are in bijective
correspondence with those of $G$ in such a way that each edge $e$ of $G$
crosses only its corresponding edge $\dual e$.
Let $\sd\fo$ be the dual spanning forest in $\dual G$ defined by $\dual e
\in \sd\fo$ iff $e \notin \fo$.
Now $G$ and $\dual G$ are unimodular since the isometry group of $\HH^2$
is unimodular, whence so are the co-compact subgroups of isometries that
fix $G$ and $\dual G$.

Note that the process of adding edges in the first
paragraph a.s.\ does not increase the set of ends of any tree: For if it
did, we could transport mass 1 from every vertex not originally in a tree
to the vertex $x$ in the tree $T$ where it joins $T$. This would give $x$
infinite mass whenever there is a new ray beginning at $x$ that uses no
edge from $T$. By the Mass-Transport Principle, the probability of this
event is 0.

Let $\fo_3$ be the set of trees in $\fo$ that have at least 3 ends.
If $\P[\fo_3 \ne \emptyset] > 0$, then by conditioning on $\fo_3 \ne
\emptyset$, we may assume (temporarily) that
this probability is 1. Now simple random walk on the
forest $\fo_3$ a.s.\ has positive speed in the metric of the forest (as we
noted above or, e.g., by Theorem 16.4 of \ref b.LP:book/), whence it also
has positive speed in the graph metric of $G$ by Lemma 4.6 of \ref
b.BLS:pert/, and hence also in the hyperbolic metric.
Therefore, as in the proof of Theorem 4.1 of \ref b.BS:hp/, simple random
walk on the forest $\fo_3$ a.s.\ converges to an ideal boundary point.
It follows that a.s.\ for every tree $T \in \fo_3$, we have that
$\mu_T$-a.e.\ ray converges to an ideal point, where
$\mu_T$ is harmonic measure on the boundary of $T$. 

Let $A$ be the limit set of $\fo$. 
We have just shown that $A \ne \emptyset$ a.s.
It follows that $A$ is dense a.s.
Indeed, for every $\epsilon > 0$, choose $B_\epsilon$ to be a non-dense
subset of the ideal boundary for which $\P[A \cap B_\epsilon \ne \emptyset]
> 1 - \epsilon$.
Now let $B$ be any non-empty open subset of the ideal boundary.
There is an isometry of $\HH^2$ that induces an automorphism of
$G$ and carries $B_\epsilon$ into $B$. Since our probability measure is
invariant under automorphisms, it follows that
$\P[A \cap B \ne \emptyset] > 1 - \epsilon$.
Since this holds for every $\epsilon > 0$, we deduce that
$\P[A \cap B \ne \emptyset] = 1$.
Since the ideal boundary is separable, the claim follows.

%if not, then let $B_1, B_2, \ldots$ be a countable set of arcs in
%the ideal boundary that generate the Borel topology.
%There must then be some $i$ such that $\P[A \cap B_i \ne \emptyset] <
%1/2$.
%This implies that $\P[A \cap B_i^{\rm c} \ne \emptyset] > 1/2$.
%However, there is an isometry of $\HH^2$ that induces an automorphism of
%$G$ and carries $B_i^{\rm c}$ into $B_i$; since our probability measure is
%invariant under automorphisms, this is a contradiction.

The density of $A$ now implies that all rays converge a.s., whence $A$ is
the entire ideal boundary.
To see this, consider any two points $\xi \ne \eta \in A$. 
Let $P, Q$ be rays that
converge to $\xi, \eta$, respectively. Let $S$ be a path between the
initial vertices of $P$ and $Q$.
No ray can cross $P \cup Q \cup S$ infinitely many times,
whence its limit set must be contained in one of the two closed
arcs determined by $\xi$ and $\eta$.
Our freedom in choosing $\xi$ and $\eta$ from the density of $A$ now gives
the result.
The same argument shows that the limit set of $\sd\fo$ is the entire ideal
boundary a.s.

If the number of trees in $\fo$ with one end is at least 3, then either
$\fo$ or $\sd\fo$
must contain a tree with at least 3 ends.
Hence, we may again conclude that the limit set of $\fo$ is the entire
ideal boundary a.s.
On the other hand, the probability is 0 that
the number of trees with one end is positive and finite, since if the
probability is positive, then as before we may assume that the entire
forest consists of such trees a.s.
In fact, by choosing just one of the trees at random, we may assume that
there is only one tree with one end.
This contradicts Theorem 5.3 of \ref b.BLPS:gip/ since $G$ is non-amenable.

The only remaining case, then, is that all trees in $\fo$ and $\sd\fo$
have 2 ends. 
We claim that this has probability 0.
For when they do, we can order the trees as the integers in the following
sense. Each tree $T$ in $\fo$ 
separates the plane into two pieces since it has two ends.
The dual of the edge boundary of $T$ lies in $\sd\fo$ and has two connected
components, each one being part of a tree in $\sd\fo$.
The same applies to each of those trees in turn, which means that on each
side of those trees, besides $T$, there is another tree in $\fo$ that
includes the dual of part of its edge boundary. Those two trees in $\fo$
are the ones next to $T$ in the integer ordering of all the trees in $\fo$.
This allows us to
define an invariant percolation with all clusters finite yet with
arbitrarily high marginal, contradicting non-amenability by Theorem 2.12 of
\ref b.BLPS:gip/.
To see this, call the unique bi-infinite simple path in a tree with 2 ends
the \dfn{trunk} of that tree.
If a vertex $x$ of the trunk is deleted and $y$ is in a finite component of
what is left of the tree (or $y = x$), then call $x$ the \dfn{trunk
attachment} of $y$.
Now given $\epsilon > 0$,
delete each tree of $\fo$ with probability $\epsilon$ independently and
in each tree that is left, delete each vertex on the trunk
with probability $\epsilon$
independently, and delete all vertices not on a trunk whose trunk
attachment was deleted.
Thus, each vertex is deleted with probability $\epsilon +
(1-\epsilon)\epsilon$.
It remains to show that the graph induced by the remaining vertices has no
infinite clusters a.s.
Number the trees by $\Z$ as indicated above, where we choose arbitrarily
which tree is numbered 0 and in which direction the integers increase.
Suppose that trees numbered $m$ and $m+n+1$ are deleted, while the $n$ trees
numbered $i$ are not deleted for $m < i < m+n+1$.
Consider a vertex $x_1$ on the trunk of tree number $m+1$. 
Then there is at least one vertex $x_2$ on the trunk of tree number $m+2$
such that for some $y_1$ whose trunk attachment is $x_1$ and some $y_2$
whose trunk attachment is $x_2$, there is an edge of $G$ between $y_1$ and
$y_2$.
Likewise, we may choose $x_3$ on the trunk of the tree number $m+3$ such
that some $y_3$ whose trunk attachment is $x_3$ is adjacent to a vertex
whose trunk attachment is $x_2$, etc.
If all vertices $x_1, \ldots, x_n$ are deleted, then there is a path in $G$
of deleted vertices stretching from tree number $m$ to tree number $m+n+1$.
The probability that $x_1, \ldots, x_n$ are all deleted is $\epsilon^n$
(recall that trees number $m+1, \ldots, m+n$ are not deleted).
We may choose infinitely many such sequences $\Seq{x_1, \ldots, x_n}$ that
are pairwise disjoint, so that the corresponding events that these
sequences are deleted are independent, each having the same probability
$\epsilon^n$.
Therefore, infinitely many such events occur a.s., and when they do, they
separate the remaining vertices between trees $m$ and $m+n+1$ into finite
components. 
Since this happens between each consecutive pair of deleted trees, all
components are finite a.s.
This completes the proof of the theorem.
%\msnote{It is also easy to define an invariant mean.}
\Qed

We can give some additional information about the limit points of the rays
in those trees with at least 3 ends. 

\procl t.3ends
Let $G$ be a one-ended graph embedded in $\HH^2$ such that a
group of isometries of $\HH^2$ acts quasi-transitively on $G$.
Let $\fo$ be an invariant forest in $G$ such that each tree has at least 3
ends.
A.s.\ for each tree in $\fo$,
the map from ends to limit points never maps more than 2 ends
to the same limit.
In addition, given that there is more than one tree in $\fo$,
a.s.\ for each tree,
the limit set is a perfect nowhere-dense set.
\endprocl

\proof
Consider a point $\xi$ of the ideal boundary that is a limit of at least 2
rays of some tree $T$ of $\fo$.
Among all the bi-infinite paths in $T$ both of whose ends converge to $\xi$,
there is only one, call it $P$, such that all others lie inside the closed
curve defined by $P \cup \{\xi\}$. Indeed, if not, it would follow that
$\xi$ is the only limit point of the ends of $\fo$. Since there is no
invariant probability measure on the ideal boundary, the event that such a
$\xi$ exists has probability 0.
Now each end with a limit $\xi$ contains a unique ray that starts at a
point in $P$. Let each vertex in such a ray send mass 1 to its starting
point in $P$. Then some points in $P$ get infinite mass,
so by the Mass-Transport Principle, this has probability 0.
This establishes that the map from ends to limit points a.s.\
never maps more than 2 ends to the same limit.
In particular, the limit set is a.s.\ infinite.

Now suppose that $\fo$ has
more than one tree a.s. Let $T$ be one of them.
If its limit set is not nowhere dense,
then it contains a maximal proper arc. Let $I$ be one such arc.
Among all the bi-infinite paths in $T$ both of whose ends converge to
points in $I$, 
there is only one, call it $P$, such that all others lie on the same side
of $P$ as $I$.
For if not, $I$ would be the entire limit set of $\fo$ and by choosing an
endpoint of $I$ at random, we would again obtain an invariant probability
measure on the ideal boundary, an impossibility.
Now each end with a limit in $I$ contains a unique ray that starts at a
point in $P$. Let each vertex in such a ray send mass 1 to its starting
point in $P$. Then some points in $P$ get infinite mass,
so by the Mass-Transport Principle, this has probability 0.
This establishes that the limit set is a.s.\ nowhere dense.

If the limit set 
is not perfect, then there is an isolated limit point, $\xi$,
and a vertex $x$ of $T$ such that three rays from $x$ that are disjoint
other than at $x$ have distinct limit points, one being $\xi$.
In fact, for each isolated limit point $\xi$, there is a unique such $x$
that is ``closest" to $\xi$ in that the ray from $x$ to $\xi$ contains no
other vertex with these properties. 
But then we can transport to $x$ mass 1 from each vertex on the ray from
$x$ to $\xi$ and so, by the Mass-Transport Principle, this has probability
0.
\Qed

In order to give examples where forests have rays that do not have positive
speed, it will be convenient to work first in the context given after the
statement of \ref t.rays-hyp-space/.
That is, we consider first forests in $\HH^d$ whose law is
invariant under isometries of $\HH^d$.
The lengths of edges are bounded and the vertices are separated
by a minimum distance in $\HH^d$.

Note that one cannot exhibit a trivial example of a forest with zero-speed
rays by subdividing edges
in a random forest in such a way that the number of subdivisions has
infinite mean, for then there is no way to re-embed the forest in an
invariant fashion. 
Our examples were discovered in conversation with David Fisher.
Take a random invariant collection of disjoint horoballs. For example,
in $\HH^2$, one can apply a random isometry to
the Ford circles in the upper half-plane model; this is possible since the
stabilizer of the set of Ford circles has co-finite measure in the full
isometry group. A fundamental domain (up to rotation in the tangent bundle)
is shown in \ref f.Ford-fund/ for the Poincar\'e disc model of $\HH^2$.
The fundamental domain has area $\pi/3$ since it is composed of two
congruent geodesic triangles, each of which has angles of measure $\pi/3$,
$\pi/2$, and $0$. 
For general $d \ge 2$, \ref b.GarRag/ show that for every discrete
subgroup $\Gamma$ of isometries of $\HH^d$ whose quotient $\HH^d/\Gamma$ is
non-compact with finite volume, there exists a $\Gamma$-invariant
collection of disjoint horoballs in $\HH^d$. For a proof that
such subgroups $\Gamma$ exist, see Chapter 14 of \ref b.Rag:book/.
By choice of $\Gamma$, we then
obtain a probability measure on collections of disjoint horoballs that is
invariant under all isometries of $\HH^d$.
A horosphere is geodesically flat, so just inside of each
horoball and at a fixed distance $\delta$ from the horosphere, we put a copy of
$\Z^{d-1}$ lying on another horosphere.
Note that we are free to choose any distance $\delta$ we want.
Again, this is done in an
invariant way. Now take any random
invariant forest IID in each copy of $\Z^{d-1}$.
None of the rays have positive speed, though all converge.
See \ref f.horoforest/ for three examples in $\HH^2$.
%The probability that a point lies inside a horodisc is $3/\pi$ since an
%easy calculation in the UHP shows that that part has area 1.
%I chose the median point on the midline.
(Such examples can be constructed similarly in complex and quaternionic
hyperbolic spaces, as well as the octonionic plane.
Although the horospheres are not then flat, one can take
them to be $\delta$-separated for any $\delta > 0$ and one may use
an embedded Cayley graph of the stabilizer of a horosphere, a group that is
finitely generated, in which one can take, say, a minimal spanning forest
independently in each horosphere.)

\efiginslabel Ford-fund {The Ford horocycle tiling with a fundamental
domain.}  x3 

\midinsert
\bigskip
      \line{%
      \hbox to 2.15truein{\hfill\Size x2 \epsfbox{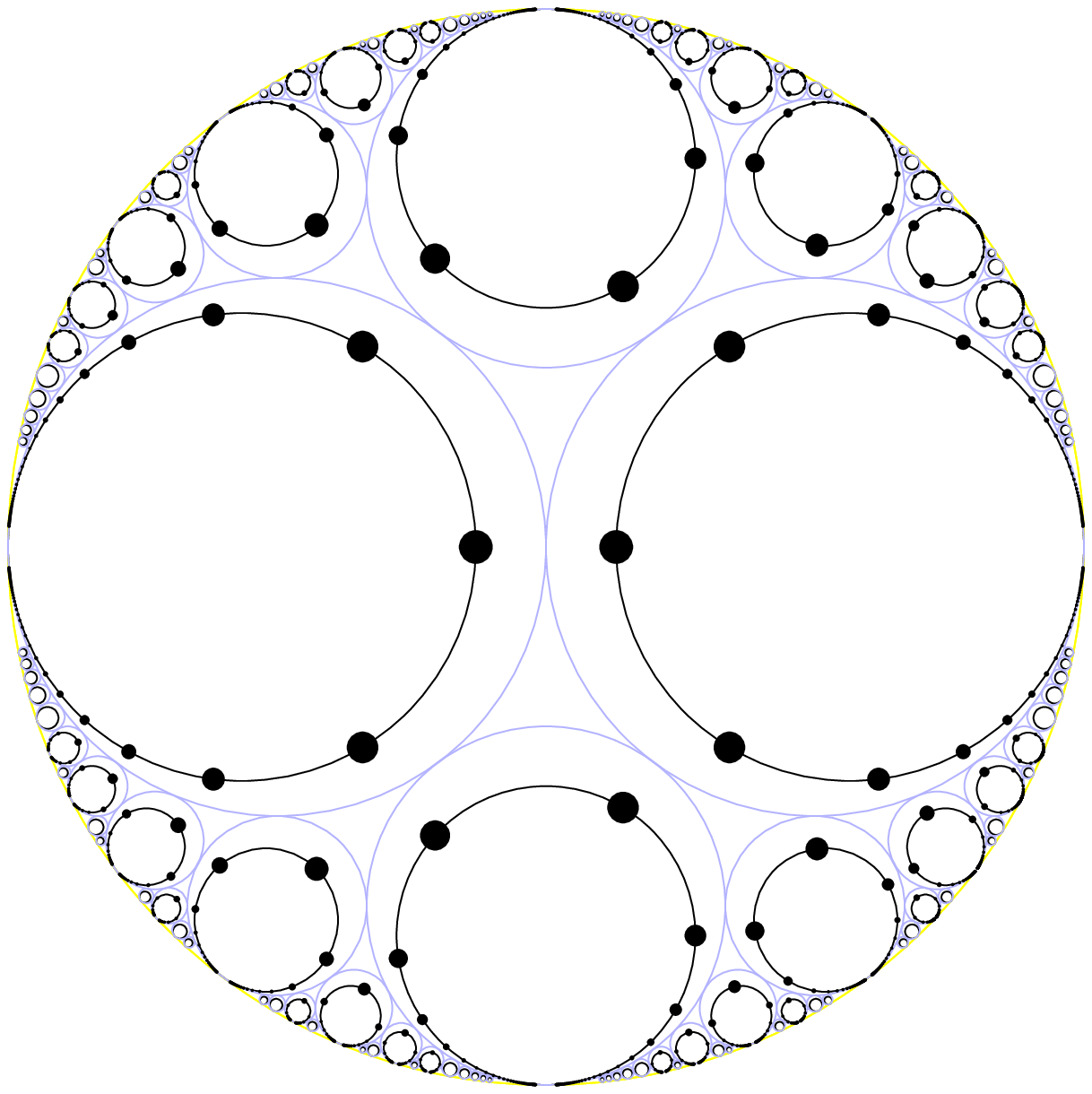}\hfill}%
      \hbox to 2.15truein{\hfill\Size x2 \epsfbox{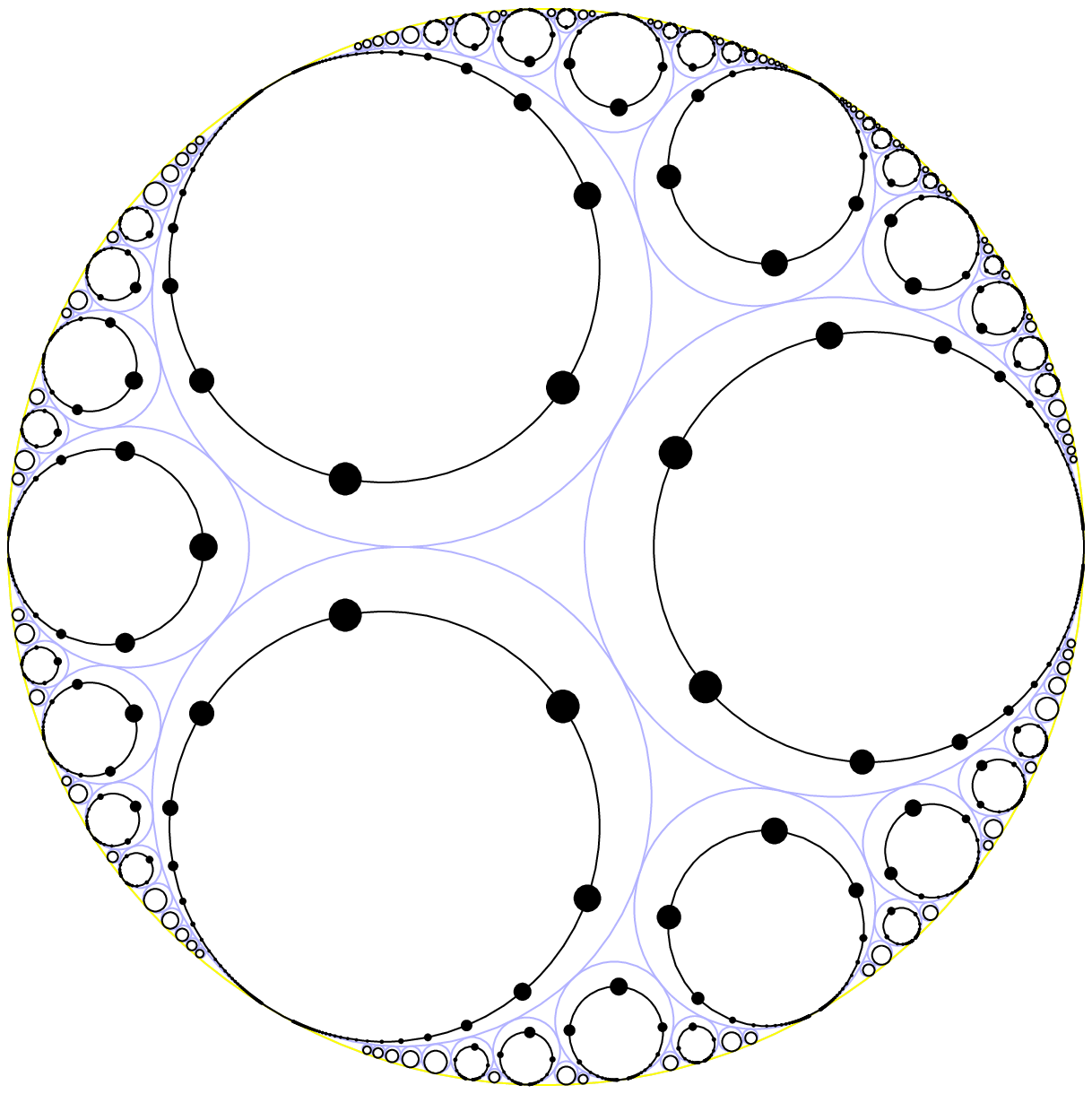}\hfill}%
      \hbox to 2.15truein{\hfill\Size x2 \epsfbox{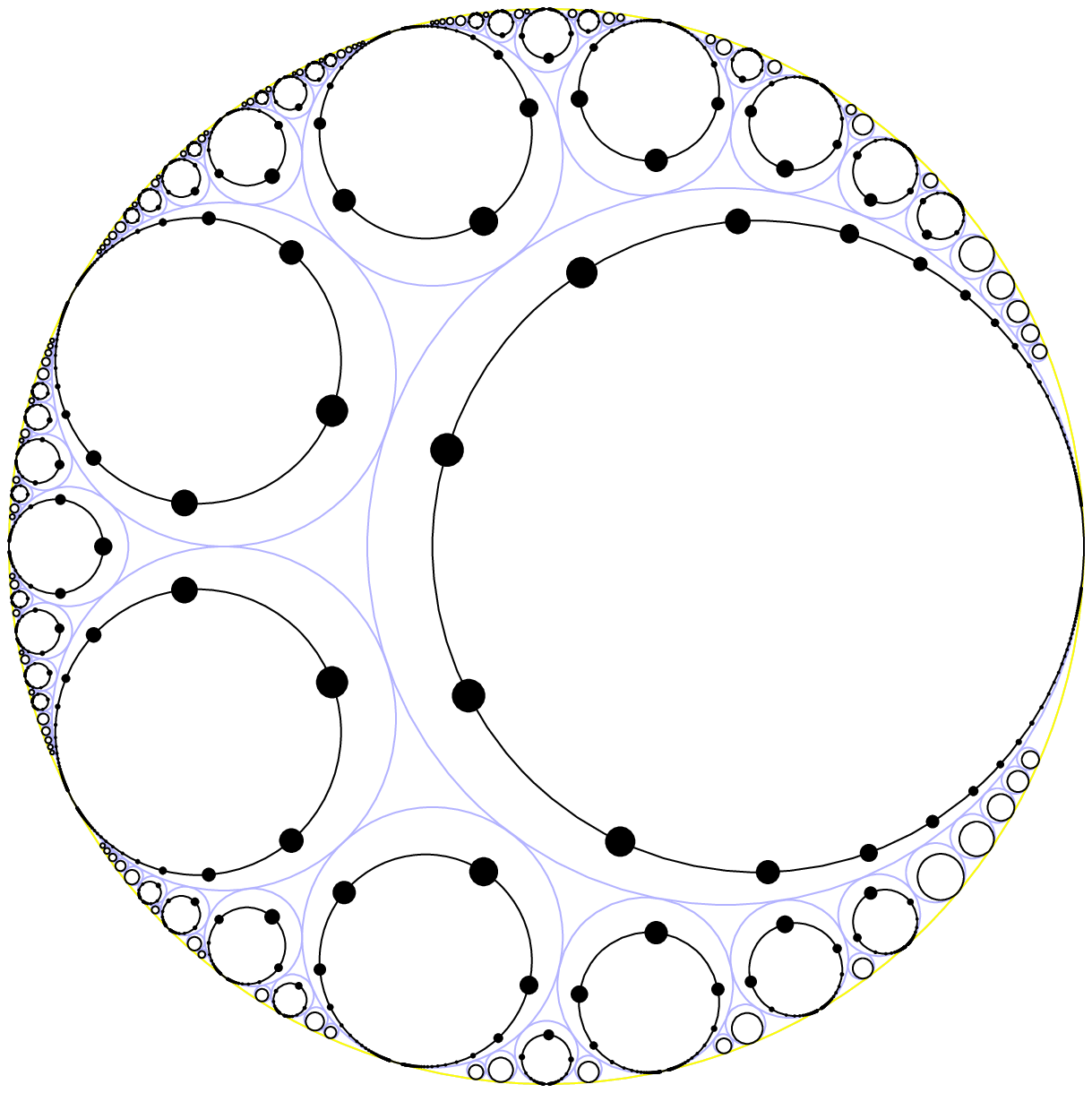}\hfill}}%
\caption{\hsize=4.9truein%
\vtop{\noindent \figlabel{horoforest}\enspace 
   The forest is in black.  The Ford horocycles are in light blue. The left
   and middle figures show examples from the corners of the fundamental
   domain, while the right figure shows a more typical example.
}}
\endinsert

Now we transfer such examples to the setting of a graph
$G$ that is quasi-isometric to $\HH^d$ and is
randomly embedded in $\HH^d$ in an isometry-invariant fashion.
For each vertex $x$ of the forest, let $\phi(x)$ be a vertex in $G$ that is
nearest to $x$; a.s.\ there is only one choice for $\phi(x)$.
If $x$ and $y$ are neighbors in the forest, then join $\phi(x)$ to
$\phi(y)$ by a shortest path in $G$; when
there is more than one choice, choose at random uniformly and
independently.
Choose the distance $\delta$ so that $\phi$ is injective and
these shortest paths are disjoint when
they come from distinct edges.
In this fashion, $\phi$ induces an embedding of the forest to a subgraph in
$G$. The choice of $\delta$ ensures that this subgraph is a forest as well.
It is the desired example.

\medskip
\noindent
{\bf Acknowledgment.} We are grateful to Omer Angel for simplifications to
some of our proofs and to 
David Fisher for permission to include the example of zero-speed rays in
hyperbolic space.

\def\noop#1{\relax}
\input \jobname.bbl

\filbreak
\begingroup
\eightpoint\sc
\parindent=0pt\baselineskip=10pt

Mathematics Department,
The Weizmann Institute of Science,
Rehovot 76100, Israel
\emailwww{Itai.Benjamini@weizmann.ac.il}
{http://www.wisdom.weizmann.ac.il/\string~itai/}

Department of Mathematics,
831 E 3rd St,
Indiana University,
Bloomington, IN 47405-7106 USA
\emailwww{rdlyons@indiana.edu}
{http://mypage.iu.edu/\string~rdlyons/}

\endgroup


\def\cprime{$'$} \def\lfhook#1{\setbox0=\hbox{#1}{\ooalign{\hidewidth
  \lower1.5ex\hbox{'}\hidewidth\crcr\unhbox0}}} \def\cprime{$'$}
  \def\cprime{$'$} \def\cprime{$'$}
\def\temp{\let\linkit=\linkyear \apaliketrue}
\temp
\ifcitationgeneration\immediate\write\labelfile{\sanitize\temp}\fi
\def\startreferences{
 \vskip0pt plus.3\vsize \penalty -150 \vskip0pt
 plus-.3\vsize \bigskip\bigskip \vskip \parskip
 \begingroup\baselineskip=12pt\frenchspacing
 \bibliographytitle
 \vskip12pt\parindent=0pt
 \def\and{{\rm and}}
 \def\em{\it}
 \def\newblock{\hskip .11em plus.33em minus.07em}
 \def\bibauthor##1{{\sc ##1}}
 \def\bibitem[##1]##2
 {\htmlanchor{##2}{}\RefLabel{##2}[##1]\hangindent=.8cm\hangafter=1}
 }
\def\endreferences{\bigskip\bigskip\endgroup}
\ifundefined{bibstylemodification}\relax\else\bibstylemodification\fi
\startreferences

\bibitem[Aizenman and Warzel (2006)]{MR2329431}
\bibauthor{Aizenman, M. \and{} Warzel, S.} (2006).
\newblock The canopy graph and level statistics for random operators on trees.
\newblock {\em Math. Phys. Anal. Geom.} {\bf 9}, 291--333 (2007).

\bibitem[Aldous and Lyons (2007)]{AL:urn}
\bibauthor{Aldous, D.J. \and{} Lyons, R.} (2007).
\newblock Processes on unimodular random networks.
\newblock {\em Electron. J. Probab.} {\bf 12}, no. 54, 1454--1508 (electronic).

\bibitem[Aldous and Steele (2004)]{MR2023650}
\bibauthor{Aldous, D.J. \and{} Steele, J.M.} (2004).
\newblock The objective method: probabilistic combinatorial optimization and
  local weak convergence.
\newblock In Kesten, H., editor, {\em Probability on Discrete Structures},
  volume 110 of {\em Encyclopaedia Math. Sci.}, pages 1--72. Springer, Berlin.
\newblock Probability Theory, 1.

\bibitem[Angel and Schramm (2003)]{MR2013797}
\bibauthor{Angel, O. \and{} Schramm, O.} (2003).
\newblock Uniform infinite planar triangulations.
\newblock {\em Comm. Math. Phys.} {\bf 241}, 191--213.

\bibitem[Barlow and Masson (2010)]{BarMass:LERW2D}
\bibauthor{Barlow, M.T. \and{} Masson, R.} (2010).
\newblock Exponential tail bounds for loop-erased random walk in two
  dimensions.
\newblock {\em Ann. Probab.} {\bf 38}, 2379--2417.

\bibitem[Benjamini and Curien (2012)]{BC:stat}
\bibauthor{Benjamini, I. \and{} Curien, N.} (2012).
\newblock Ergodic theory on stationary random graphs.
\newblock {\em Electron. J. Probab.} {\bf 17}, Article 93, 20 pp.\
  (electronic).

\bibitem[Benjamini, Lyons, Peres, and Schramm (1999)]{MR99m:60149}
\bibauthor{Benjamini, I., Lyons, R., Peres, Y., \and{} Schramm, O.} (1999).
\newblock Group-invariant percolation on graphs.
\newblock {\em Geom. Funct. Anal.} {\bf 9}, 29--66.

\bibitem[Benjamini, Lyons, and Schramm (1999)]{BLS:pert}
\bibauthor{Benjamini, I., Lyons, R., \and{} Schramm, O.} (1999).
\newblock Percolation perturbations in potential theory and random walks.
\newblock In Picardello, M. \and{} Woess, W., editors, {\em Random Walks and
  Discrete Potential Theory}, Sympos. Math., pages 56--84, Cambridge. Cambridge
  Univ. Press.
\newblock Papers from the workshop held in Cortona, 1997.

\bibitem[Benjamini and Schramm (2001a)]{BS:hp}
\bibauthor{Benjamini, I. \and{} Schramm, O.} (2001a).
\newblock Percolation in the hyperbolic plane.
\newblock {\em Jour. Amer. Math. Soc.} {\bf 14}, 487--507.

\bibitem[Benjamini and Schramm (2001b)]{MR1873300}
\bibauthor{Benjamini, I. \and{} Schramm, O.} (2001b).
\newblock Recurrence of distributional limits of finite planar graphs.
\newblock {\em Electron. J. Probab.} {\bf 6}, no. 23, 13 pp. (electronic).

\bibitem[Bollob{\'a}s and Riordan (2011)]{MR2839983}
\bibauthor{Bollob{\'a}s, B. \and{} Riordan, O.} (2011).
\newblock Sparse graphs: metrics and random models.
\newblock {\em Random Structures Algorithms} {\bf 39}, 1--38.

\bibitem[Bowen (2003)]{MR2026846}
\bibauthor{Bowen, L.} (2003).
\newblock Periodicity and circle packings of the hyperbolic plane.
\newblock {\em Geom. Dedicata} {\bf 102}, 213--236.

\bibitem[Connes, Feldman, and Weiss (1981)]{MR84h:46090}
\bibauthor{Connes, A., Feldman, J., \and{} Weiss, B.} (1981).
\newblock An amenable equivalence relation is generated by a single
  transformation.
\newblock {\em Ergodic Theory Dynamical Systems} {\bf 1}, 431--450 (1982).

\bibitem[Elek (2010)]{MR2776719}
\bibauthor{Elek, G.} (2010).
\newblock On the limit of large girth graph sequences.
\newblock {\em Combinatorica} {\bf 30}, 553--563.

\bibitem[Elek and Lippner (2010)]{MR2566316}
\bibauthor{Elek, G. \and{} Lippner, G.} (2010).
\newblock Sofic equivalence relations.
\newblock {\em J. Funct. Anal.} {\bf 258}, 1692--1708.

\bibitem[Gaboriau and Lyons (2009)]{MR2534099}
\bibauthor{Gaboriau, D. \and{} Lyons, R.} (2009).
\newblock A measurable-group-theoretic solution to von {N}eumann's problem.
\newblock {\em Invent. Math.} {\bf 177}, 533--540.

\bibitem[Garland and Raghunathan (1970)]{MR0267041}
\bibauthor{Garland, H. \and{} Raghunathan, M.S.} (1970).
\newblock Fundamental domains for lattices in ({R}-)rank {$1$} semisimple {L}ie
  groups.
\newblock {\em Ann. of Math. (2)} {\bf 92}, 279--326.

\bibitem[Hjorth (2006)]{MR2258624}
\bibauthor{Hjorth, G.} (2006).
\newblock A lemma for cost attained.
\newblock {\em Ann. Pure Appl. Logic} {\bf 143}, 87--102.

\bibitem[Kapovich and Benakli (2002)]{MR1921706}
\bibauthor{Kapovich, I. \and{} Benakli, N.} (2002).
\newblock Boundaries of hyperbolic groups.
\newblock In Cleary, S., Gilman, R., Myasnikov, A.G., \and{} Shpilrain, V.,
  editors, {\em Combinatorial and Geometric Group Theory}, volume 296 of {\em
  Contemp. Math.}, pages 39--93. Amer. Math. Soc., Providence, RI.
\newblock Papers from the AMS Special Sessions on Combinatorial Group Theory
  and on Computational Group Theory held in New York, November 4--5, 2000 and
  in Hoboken, NJ, April 28--29, 2001.

\bibitem[Lov\'asz (2012)]{Lovasz:limits}
\bibauthor{Lov\'asz, L.} (2012).
\newblock {\em Large Networks and Graph Limits}.
\newblock American Mathematical Society, Providence, RI.


\bibitem[Lyons with Peres (2013)]{LP:book}
\bibauthor{Lyons, R. {\rm with} Peres, Y.} (2013).
\newblock {\em Probability on Trees and Networks}.
\newblock Cambridge University Press.
\newblock In preparation. Current
  version available at 
  %\hfil\break
  %\htmlref{http://mypage.iu.edu/\string~rdlyons/}{{\tt
  %http://mypage.iu.edu/\string~rdlyons/}}.
  \url{http://mypage.iu.edu/~rdlyons/}.

\bibitem[Pemantle (1991)]{MR92g:60014}
\bibauthor{Pemantle, R.} (1991).
\newblock Choosing a spanning tree for the integer lattice uniformly.
\newblock {\em Ann. Probab.} {\bf 19}, 1559--1574.

\bibitem[Raghunathan (1972)]{MR0507234}
\bibauthor{Raghunathan, M.S.} (1972).
\newblock {\em Discrete Subgroups of {L}ie Groups}.
\newblock Springer-Verlag, New York.
\newblock Ergebnisse der Mathematik und ihrer Grenzgebiete, Band 68.

\bibitem[Teplyaev (1998)]{MR1658094}
\bibauthor{Teplyaev, A.} (1998).
\newblock Spectral analysis on infinite {S}ierpi\'nski gaskets.
\newblock {\em J. Funct. Anal.} {\bf 159}, 537--567.

\endreferences
\bye